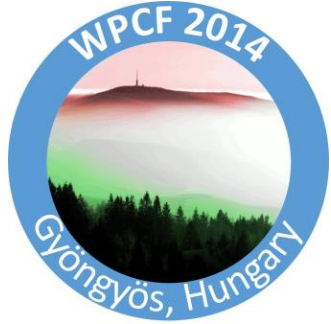

# Educational Challenges of Rubik's Cube

Sándor Kiss[1]

[1]Rubik Studio, Városmajor 74., Budapest,
Hungary, H-1122

March 31, 2015

**Abstract**

The first 2x2x2 twisty cube was created as a demonstration tool by Ernő Rubik in 1974 to help his students understand the complexity of space and the movements in 3D. He fabricated a novel 3x3x3 mechanism where the 26 cubies were turning, and twisting independently, without falling apart. The cube was dressed in sophisticated colors which made it a unique puzzle. Even without instruction is the aim of the game was self-explanatory. Its educational value in VSI (Visual-Spatial Intelligence), developing strategy, memorization and logistics, improve concentration and persistence in problem solving is high in every age group. A logical puzzle has outreach far beyond. Those aspects are briefly covered in this article.

## 1 Introduction

The role of experience in learning has been researched since long. Let us pick up education, mathematics and puzzle. The relations among those 3 factors will or may determine the fate of a puzzle. The relations namely puzzle and education; mathematics and education; mathematics and puzzle; have been showing different behaviors on various puzzle objects. A well balanced mathematical and educational harmony of the puzzle is necessary, but not enough. The puzzle needs brain work and concentration, should give amusement, it should also engage the player and its attractiveness is preferable. In addition the puzzle's usability as educational tool is substantial value. Whether a new puzzle will bring success or will fail the answer has been occultly determined above relations. The difficult point is how to interpret the message of puzzle-education-mathematics triad for a given puzzle?

Let us talk about the ordinary cube! Even a baby meets with the cube in early age and likes to play with it, without knowing that it is a solid with six congruent square faces, a regular hexahedron which has edges, faces, vertex, face diagonals and space diagonals, just no to go further except one more characteristic the symmetry, the beauty of cube. For a baby, for a young child the cube is fun, a challenging interesting shape. Two or more cubes are more amusing, it gives a rich feeling. Each cube can be personalized by dressing, can be used for construction, and can be moved, turned and thrown. As the toddlers are coming into kindergarten and school age they meet more cubes in various forms like a lump of sugar, a cake, cube cut food, a dice, a cube puzzle for kids, a furniture, a rock salt crystal, plus boxes, frames, ornaments and many cube shaped equipment and utensils. The cubes before everything else are good for toys and perfect for mind teasers. When children go to school and learning art, geometry and mathematics they may meet with variations of cube puzzles: soma, jigsaw cube, secret cubes, folding cubes and during the last four decades also the Rubik's Cube, the originator of family of twisty puzzles.

We can say in the kindergarten the cube is a toy, in elementary school the cube is more an object to study and its playfulness having been pushed into the background. Moreover in junior high schools the math of the cube receives priority. In this learning process familiarity is key factor for the teacher, and this notoriety greatly helps in transmitting the abstract knowledge of cube. This is the way how we did and they acquire the algebra and geometry of the cube, the Platonic object. The teaching comes with learning cube language as well i.e.; 6 cubed is 216, cube root of 27 is 3, the perfect cubes and so on, some day or other follows with more advanced descriptions e.g. octahedron is the dual polyhedron to the cube. However symmetry brings us to highest level, to isometries and group theories of the cube. But the main point is the better you know the more easily you learn. The teacher can also help the deeper understanding by offering choice between scientific approach and heuristic techniques. The steps upward in learning are student, teacher and science.

## 2 History

After having said that let us talk about Rubik's Cube. It is well known that the first 2x2x2 twisty cube was created as a demonstration tool by Ernő Rubik in 1974 to help his students understand the complexity of space and the movements in 3D. By turns of individual, but interrelated cubes changed their positions and a magical scrambling occurred. A new object was born, but no strategy existed how to put back the little cubes into starting order using. Ernő was both emotionally and intellectually captured by it. There was a structural problem, because the elastic band to keep the cubes together produced a knot and was torn soon. It was the next challenge for Ernő to fabricate a mechanism to overcome that problem. His novel construction gave birth to 3x3x3 cube. Now the cube layers were turning smoothly. The 26 cubies were turning, and twisting independently, without falling apart. Its symmetry implies perfectness.[1] A genuine invention, that was an object supposed to be impossible before. The patent for his Spatial Logical Toy was granted within a short time on 28th of October, 1976. The cube was dressed in sophisticated colors which turned it into unique puzzle. Even without instruction is the aim of the game was self-explanatory.

The cube was ready, but how to put back the colors was unknown. Here was the next provoking challenge! Experimenting, perhaps better to say a discovery was necessary. Ernő had started the journey on unknown road and after had been into quite a few standoffs and dead ends he has arrived to a safe solution several weeks later. His strategy and his set of algorithms were necessary and enough to solve the cube from any scrambled position.

Of course the potential of his cube had been sensible more and more and gave vigorous reason to get it manufactured as soon as possible. But the cube had been born in a small country existing under communist regime. We "enjoyed" all the drawbacks of state monopolized companies and foreign trade behind the Iron Curtain. However Politechnika, a toy-industrial cooperative after nine months consideration was brave enough to sign a contract with Ernő. Production started and 12,000 pieces *Magic Cube* (Bűvös kocka) were produced in 1977. The official trading companies were not much interested; according to their opinion it was unsolvable, unsalable, unusual piece of junk. The toyshops cautiously placed order for 5000 cubes and without big noise all were sold within short time.

It happened like in the Aladdin fairy tail the lamp had rubbed inadvertently the Genie of the cube has been unleashed. Soon it was in the schools, also seized by teachers during lessons, and as present was given to scientists at conferences outside of Hungary. The Magic Cubes reached the world. The international commercial attempts often ended in vain, but gave a demanding challenge to do it! The scientists loved it and also faced such challenges: What kind of mathematics is hidden in it? Is there any similar function in physics, chemistry or nature? The support and enthusiasm from scientists was important, but it was just observed among buyers but hardly convinced them. Finding a potential business partner was a see-saw situation where one side changed frequently, but in 2 years Ideal Toys, USA went for it. The product as *Rubik's Cube* reached its media boiling point at New York Toy Fair in 1980 and right away took the world by storm.

## 3 Solution

Those were the exiting days when everyone had to find out the solution by oneself. The solution by manual era arrived next which was followed by the current web tutorial era. Scientifically the first step was the notation of the cube developed by David Singmaster to denote the sequence of moves. The diffusion of the cube among scientists was quick and raised theoretical excitement. Mathematical models and computer programs were developed; analogies in physics and among atomic particles just to name a few were pointed out. The appearance of cube overlapped the peak period in group theory research. Soon Rubik's Cube was selected as illustration of permutation groups. Almost simultaneously the society was shaken by wave of Rubik and cube stories, jokes, songs, cartoons, films, TV series and many more off sides. But definitely it also reached education institutes quite early. Nowadays on the web anyone can select from the plethora of solutions online, surf the social impacts of the cube, wonder the effects on art, architecture, literature, fashion and competitions, follow the records connected to Rubik's Cube and enjoy scientists' talks about it. A huge collection is waiting for correct categorization and systematic processing.

## 4  In Education

Scientists who were dealing with Rubik's Cube gave lectures and composed tasks with it for exercise. That was a timely initiative as most of the students were perplexed already by cube. At university level of education each explanation generated request to answer more new queries. It became a compulsory topic in group theory and abstract algebra curriculums. Search for the minimum necessary moves was on high before long. In 2010 utilizing idle time of Google's servers was concluded that God's Number is 20, but it was done by programmed brutal force but the proof is still missing. [2]

Let us proceed to middle and elementary school. Teachers today attracted to the puzzle-solving lesson because it helps to learn mathematics and develops direction-following and memorization plus persistence. But there was a complete different standpoint at the beginning. The Hungarian National Institute of Pedagogy turned down Ernő's initiative to use it in education in 1976 as improper teaching tool. After having the great 1980 success more and more teacher set out to use with good results, but the sales of the cube dropped to negligible level at mid eighties and the speedcubing competitions give a little bust and kept it alive during the nineties. This also downward affected the cube's role in education. Also the torrent of new twisty puzzles drove teacher's intentions into a corner; latest ones might be more attractive than good old cube.

By turn of the century speedcubing was gaining popularity all over the world due to new teenage generation with strong competition spirit supported by wide use of internet. Also it was a decisive factor that the World Cube Association had been established in 1999. WCA laid down standards and rules, governs competitions and keeps records for all puzzles labeled as Rubik.

After celebrating the 30th birthday of Rubik's Cube it became opportune to do something for education. During the past three decades teaching methods have developed and experiments in schools here and there have accumulated a lot of results. Furthermore we have seen new cube like puzzles, for example the void cube without center cubies on all three axes, or giant cubes with as many as 6, 7 or even 10 cubies across cube faces since 2000. Also the family of odd shapes twisty puzzles was grown by dodecahedron, bandaged cubes and so on. Many interesting puzzles, but they have received less and less attention on the market.

Ernő always supported the introduction of the cube to education at all levels. The spatial awareness is still not strong enough in schools worldwide. Most of the spatial teaching tools to meet are of classical unmoving type. It makes a big difference seeing 3D demos versus taking a workable model into own hands. The latter one builds confidence and also helps underachieving students. His initiative and long wish was taken up by Seven Towns, the company behind the world success of the Rubik's Cube. After several year of preparation the project You CAN do the Rubik's Cube program was launched in 2010. This complex program is to help students across the U.S. and beyond realize their potential while emphasizing teamwork, "outside of the box" thinking, and creativity.[3] Now the Rubik's Cube is in the STEAM (Science, Technology, Engineering, Art, Mathematics) project, which offers the teachers MATH methodology: M – Modify the lesson for understanding; A – Apply the lesson to repeated practice; T– re-Teach the lesson for mastery; H – Higher level learning for enrichment. The need for You CAN do the Rubik's Cube project is justified by 31,000+ implemented packages so far.

## 5   Summary

The 40th birthday brought the internationally-acclaimed Beyond Rubik's Cube exhibition at Liberty Science Center in New Jersey in 2014 and also we said Happy Returns of the Day to Ernő on his 70th birthday on 13th of July. The BRC exhibition will travel to several cities around the world (Cleveland and Edmonton, Canada in 2105) and hopefully after 5-7 years it will arrive to Budapest. Its interactive content will be monitored, changed and renewed as changes of time would request. We can just repeat Julius Ceasar's famous phrase Alea iacta est, i.e. events have passed a point of no return. The Rubik's Cube now has been living his own life.

One can hardly close this ever developing story, because new developments have been coming since 1974. Instead of closing remark let me quote John von Neumann:

> "By and large it is uniformly true that in mathematics that there is a time lapse between a mathematical discovery and the moment it becomes useful; and that this lapse can be anything from 30 to 100 years, in some cases even more; and that the whole system seems to function without any direction, without any reference to usefulness, and without any desire to do things which are useful." [4]

## References

[1] Marcus du Sautoy: *Symmetry, reality's riddle*. TEDGlobal lecture 2009.

[2] Tomas Rokicki, Herbert Kociemba, Morley Davidson, John Dethridge: *God's Number is 20* (2010).

[3] Rohrig, Brian: *Puzzling Science: Using the Rubik's Cube to Teach Problem Solving*. Science Teacher, Volume 77, No. 9 p54-56, Dec 2010. See Rohrig's article.

[4] John von Neumann: Collected Works vol. VI, p. 489.